\documentclass[referee]{raa}            

 \usepackage{natbib}
\usepackage{amssymb}
\usepackage{amsmath}
\usepackage{graphicx,times}             
\input{epsf.sty}                        
\input{psfig.sty}                       

\newcommand{\amc}{{(1-\mu) }}

\newcommand{\zabs}{{\Bigl(1+\frac{5A_2}{2{r^2_2}_i}\Bigr)}}

\begin{document}

   \title{Poynting-Robertson Effect on the Lyapunov Stability of Equilibrium Points in the Generalized Photogravitational Chermnykh's Problem
}

   \volnopage{Vol.0 (200x) No.0, 000--000}      
   \setcounter{page}{1}          

   \author{Badam Singh Kushvah \inst{1}}

   \institute{Department of    Mathematics,\\
National Institute of Technology,  Raipur (C.G.) -492010,INDIA\\\email{bskush@gmail.com}}

   \date{Received~~2009 month day; accepted~~2009~~month day}

\abstract{ The  Poynting-Robertson(P-R) effect on Lyapunov stability of equilibrium points is being discussed in the the generalized photogravitational Chermnykh's problem when bigger primary is a sours of radiation and smaller primary is an oblate spheriod. We  derived the equations of motion, obtained the equilibrium points and examined the linear stability of the equilibrium points for various values of parameter which have been used in the present problem. We have examined the effect of  gravitational potential from the belt. The positions of the equilibrium points are  different  from the position in classical  case. We have seen  that due to the P-R  effect all the equilibrium points are unstable in Lyapunov sense. 
\keywords{equilibrium points: generalized
photogravitational:Chermnykh's problem:radiation pressure: Poynting-Robertson effect}
}

   \authorrunning{Badam Singh Kushvah }            
   \titlerunning{Poynting-Robertson Effect  ---- Chermnykh's Problem}  

   \maketitle

%
%
\section{Introduction}
\label{intro}
The Chermnykh's problem is new kind of restricted three body problem which was first time studied by \cite{Chermnykh1987VeLen}. \cite{PapadakisKanavos2007Ap&SS}  given numerical exploration of Chermnykh's problem, in which the equilibrium points and zero velocity curves studied numerically also the non-linear stability for the triangular Lagrangian points are computed numerically for the Earth-Moon and Sun-Jupiter mass distribution when the angular velocity varies. The solar radiation pressure force $F_p$ is exactly apposite to the gravitational attraction force $F_g$ and change with the distance by the same law it is possible to consider that the result of action of this force will lead to reducing the effective mass of the  Sun or particle. It is acceptable to speak about a reduced mass of the particle as the effect of reducing its mass depends on the properties of the particle itself. The  first order in
${\vec{\frac{V}{c}}}$ the radiation pressure force is given by[see \cite{Poynting1903} \cite{Robertson1937}]:
\begin{equation}\vec{F}=\displaystyle{{F_p}\biggl\{ \frac{\vec{R}}{R} - \frac{\vec{V}.\vec{R}\vec{R}}{cR^{2}}- \frac{\vec{V}}{c} \biggr\}}\label{eq:F}\end{equation}
Where $F_{p}$=$\frac{3Lm}{16{\pi}R^{2}\rho{sc}}$   denotes the
measure of the radiation pressure force, $\vec{R}$  the position
vector of $P$ with respect to radiation sours Sun  $S$, $\vec{V}$ the
corresponding velocity vector and $c$ the velocity of light. In the
expression of $F_{p}$, $L$ is luminosity of the radiating body,
while $m$, $\rho $ and $s$ are the mass, density and cross section
of the particle respectively.

The first term in equation (~\ref{eq:F}) expresses the radiation pressure. The second term represents the Doppler shift of the incident radiation and the third term is due to the absorption and subsequent re-emission of the incident radiation. These last two terms taken together are the Poynting-Robertson effect. The Poynting-Robertson effect will operate to sweep small particles of the solar system into the Sun at cosmically rapid rate. \cite{Chernikov1970}   discussed the
position as well as the stability of the Lagrangian equilibrium
points when radiation pressure, P-R drag force are included.
~\cite{Murray1994} systematically discussed the dynamical effect of
general drag in the planar circular restricted three body problem.
~\cite{KushvahBR2006} examined the linear stability of triangular
equilibrium points in the generalized photogravitational restricted
three body problem with Poynting-Robertson drag, $L_4$ and $L_5$
points became unstable due to P-R drag which is very remarkable and
important, where as they are linearly stable in classical problem when $0<\mu<\mu_{Routh}=0.0385201$. Further the normalizations  of Hamiltonian and nonlinear stability of $L_{4(5)}$ in the present of P-R drag has been studied by \citet*{Kushvah2007BASI,Kushvah2007Ap&SS.312,Kushvah2007EM&P..101} 

In this paper we  generalized  our previous paper \cite{Kushvah2008Ap&SS} in which  linear stability has been examined in the generalized photogravitaional Chermnykh's problem, we have seen that the collinear points are linearly unstable and triangular points are stable in the sense of Lyapunov stability provided $\mu<\mu_ {Routh} =0.0385201$.  In present paper we have obtained the equations of motion and equilibrium points, zero velocity curves and the linear stability in the generalized photogravitaional Chermnykh's problem. We have found that the collinear points deviate from the axis joining the two primaries, while the triangular points are not symmetrical due to Poynting-Robertson effect. We have examined the effect of  gravitational potential from the belt, oblateness effect and radiation effect on the  Lyapunov stability.
\section{Equations of Motion and Position of Equilibrium Points}
\begin{figure}
    \plotone{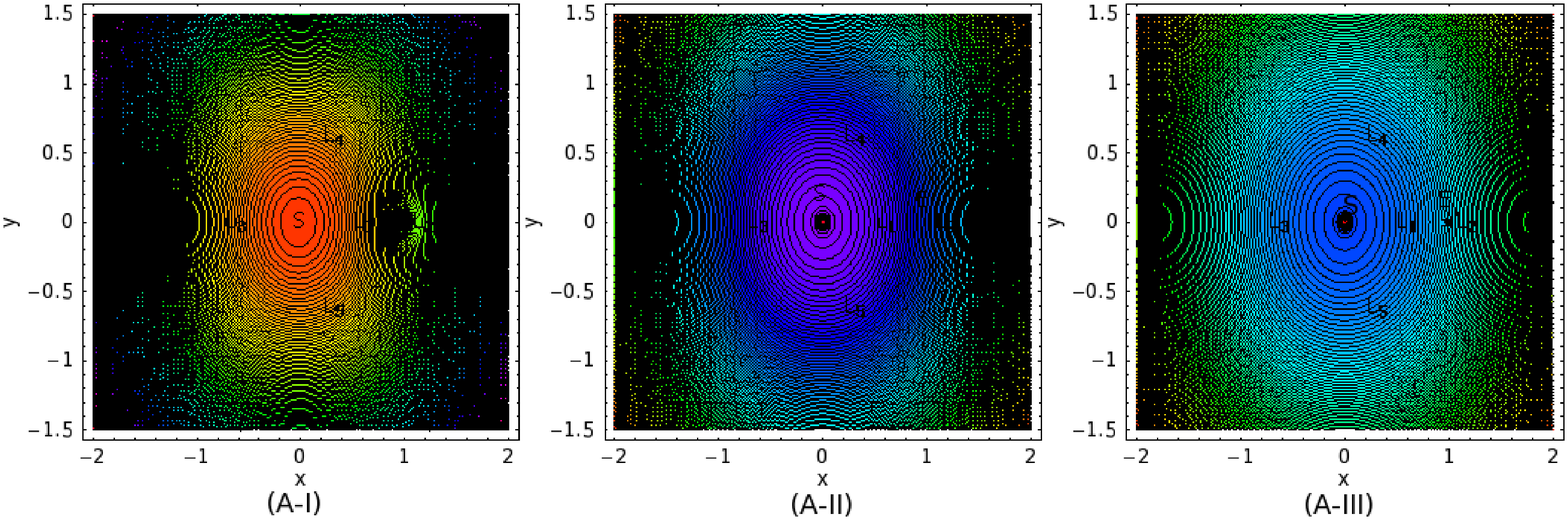} \plotone{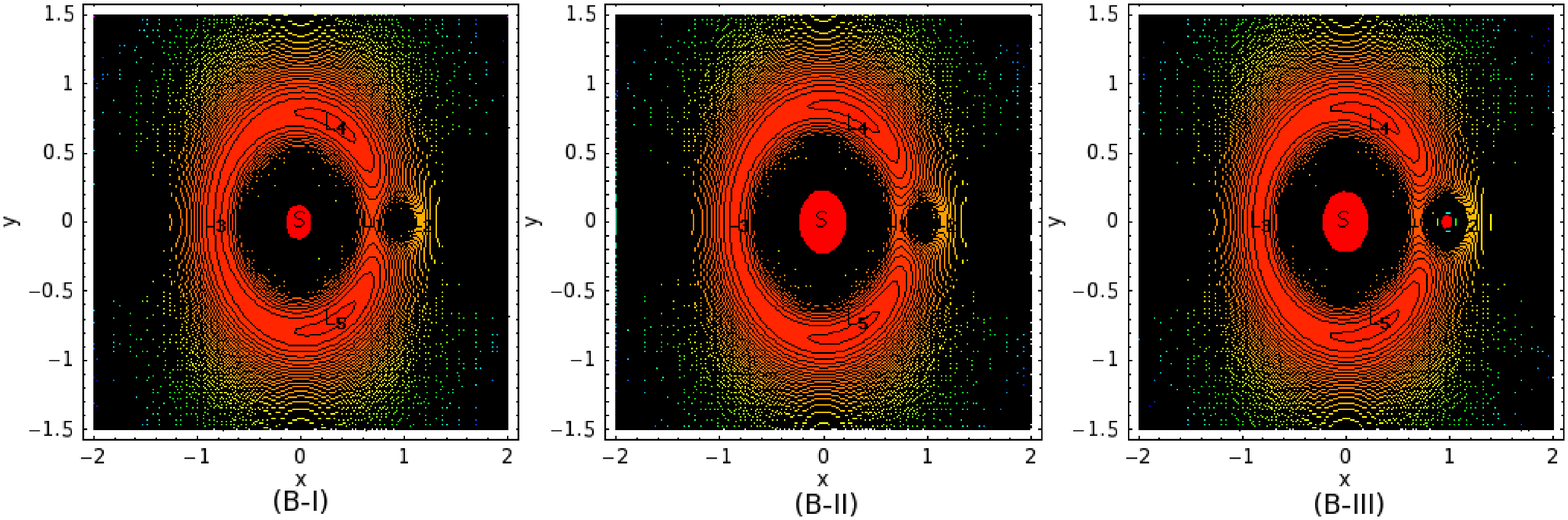} \plotone{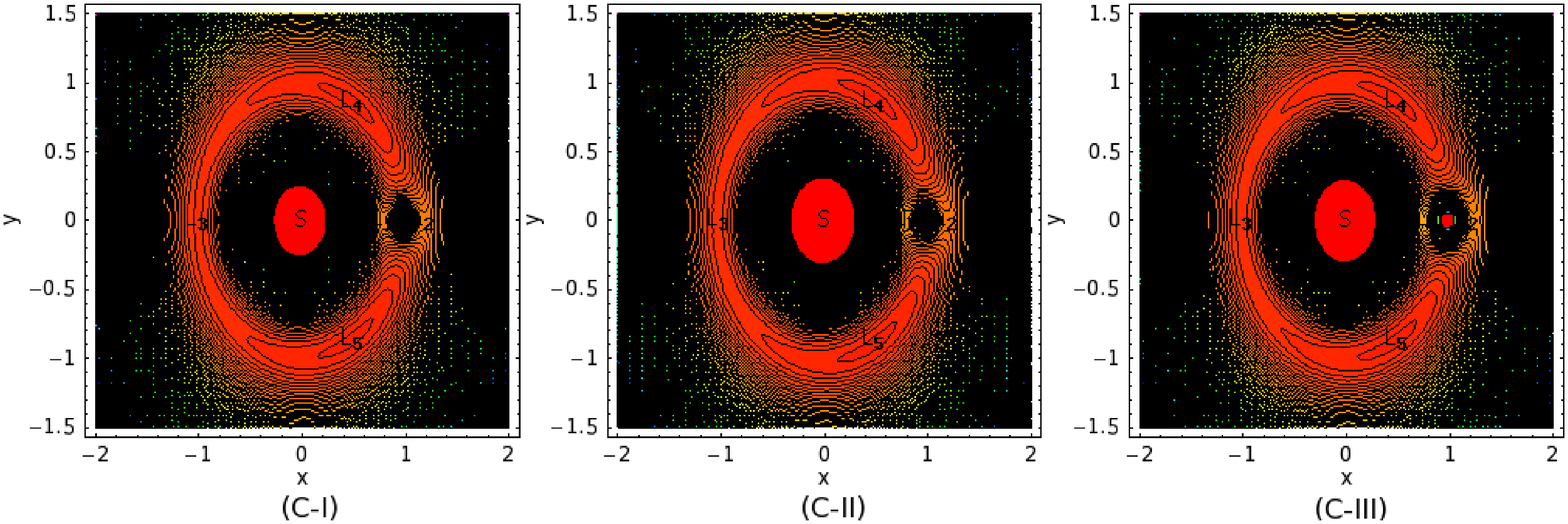} \plotone{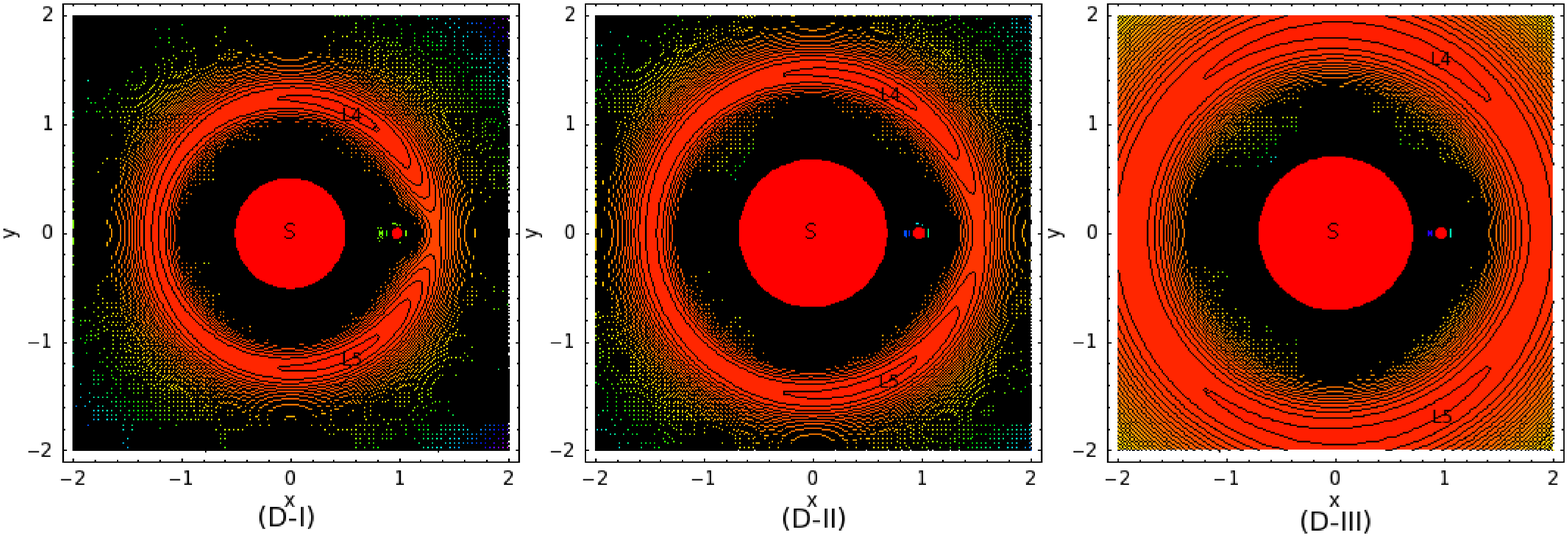}
   \caption{This figure show the zero velocity curves for $\mu=0.025, T=0.01$, frames (A-I to III)$q_1=0.00, M_b=0.02, A_2=0.00--0.04$, (B-I to III)$q_1=0.50, M_b=0.02, A_2=0.00--0.04$,(C-I to III)$q_1=1.00, M_b=0.02, A_2=0.00--0.04$,(D-I to III)$q_1=1.00, M_b=0.25, 0.5, 0.75, A_2=0.02$}
   \label{fig:plotzvc}
\end{figure}

Let us consider the model proposed by \cite{MiyamotoNagai1975PASJ}, according to this model the potential of belt is given by:
\begin{eqnarray}
 V (r,z)=\frac{\mathbf{b}^2M_b\left[\mathbf{a}r^2 +\left(\mathbf{a}+3N\right)\right]\left( \mathbf{a}+N\right)^2}{N^3\left[r^2+\left(\mathbf{a}+N\right)^2\right]^{5/2}}
\end{eqnarray}
where $M_b$ is the total mass of the belt and $r^2=x^2+y^2$, $\mathbf{a,b}$ are parameters which determine the density profile of the belt, if  $\mathbf{a}=\mathbf{b}=0$ then the  potential equals to the one by a point mass. The parameter $\mathbf{a}$  \lq\lq flatness parameter\rq\rq and $\mathbf{b}$  \lq\lq core parameter\rq\rq.  
 where $N=\sqrt{z^2+\mathbf{b}^2}$,\ $T=\mathbf{a}+\mathbf{b}$, $z=0$.
Then we obtained \begin{equation}
 V(r,0)=-\frac{M_b}{\sqrt{r^2+T^2}}\label{eq:Vr0}
\end{equation}
and  \(V_x=\frac{M_b x}{\left(r^2+T^2\right)^{3/2}},\ V_y=\frac{M_b y}{\left(r^2+T^2\right)^{3/2}} \).
As in   \cite{Kushvah2008Ap&SS, Kushvah2008Ap&SS.315}, we consider the barycentric rotating co-ordinate system $Oxyz$ relative to inertial system with angular velocity $\omega$ and common $z$--axis.  We have taken line joining the primaries as $x$--axis. Let $m_1, m_2$ be the masses of bigger primary(Sun)  and smaller primary(Earth) respectively. Let  $Ox$, $Oy$  in the equatorial plane of smaller primary  and $Oz$ coinciding with the polar axis of $m_{2}$. Let $r_{e}$, $r_{p}$ be the equatorial and polar radii of $m_{2}$ respectively,  $r$ be the distance between primaries.  Let infinitesimal mass $m$ be placed at the  point $P(x,y,0)$. We take units such that sum of the masses and distance between primaries as  unity, the unit of time i.e. time period of $m_{1}$ about $m_{2}$  consists of $2\pi$ units such that the Gaussian constant of gravitational $\Bbbk^{2}=1$. Then perturbed mean motion $n$ of the primaries is given by  $n^{2}=1+\frac{3A_{2}}{2}+\frac{2M_b r_c}{\left(r_c^2+T^2\right)^{3/2}}$, where $r_c^2=(1-\mu)q_1^{2/3}+\mu^2$. For simplicity,  we set $r=r_c=0.9999, T=0.01$ for further numerical results, where $A_{2}=\frac{r^{2}_{e}-r^{2}_{p}}{5r^{2}}$ is oblateness coefficient of $m_{2}$.
where $\mu=\frac{m_{2}}{m_{1}+m_{2}}$  is mass parameter,  $1-\mu=\frac{m_{1}}{m_{1}+m_{2}}$ with $m_{1}>m_{2}$. Then coordinates of $m_1$ and $m_2$ are  $(x_1,0)=(-\mu,0)$ and  $(x_2,0)=(1-\mu,0)$ respectively. Further, in our consideration, the velocity of light needs to be dimensionless, too, so consider  the dimensionless velocity of light as $c_d=c$ which depends on the physical masses of the two primaries and the distance between them. The mass of Sun $m_1\approx 1.989\times 10^{30}kg$ $\approx 332,946m_2$(The mass of Earth), hence mass parameter for this system is $\mu=3.00348\times10^{-6}$.  In the above mentioned reference system  the  we determined  the equations of motion of the infinitesimal mass particle in $x y$-plane.  Now using  \cite{MiyamotoNagai1975PASJ} profile and  \cite{Kushvah2008Ap&SS, Kushvah2008Ap&SS.315}, then the equations of motion are given  by:
\begin{eqnarray}
\ddot{x}-2n\dot{y}&=&U_{x}-V_x=\Omega_x ,\label{eq:Omegax}\\
\ddot{y}+2n\dot{x}&=&U_{y}-V_y=\Omega_y\label{eq:Omegay}
 \end{eqnarray}
where
\begin{eqnarray*}
&&\Omega_x= n^{2}x-\frac{(1-\mu)q_1(x+\mu)}{r^3_1}-\frac{\mu(x+\mu-1)}{r^3_2}\\&&-\frac{3}{2}\frac{\mu{A_2}(x+\mu-1)}{r^5_2}-\frac{M_b x}{\left(r^2+T^2\right)^{3/2}} \nonumber\\
&&-\frac{W_1}{r^2_1}\biggl\{\frac{(x+\mu)}{r^2_1}[(x+\mu){\dot{x}+y\dot{y}}] +\dot{x}-ny \biggr\},\\
&&\Omega_y=n^{2}y
-\frac{(1-\mu)q_{1}{y}}{r^3_1}
-\frac{\mu{y}}{r^3_2}\nonumber\\&&-\frac{3}{2}\frac{\mu{A_2}y}{r^5_2}-\frac{M_b y}{\left(r^2+T^2\right)^{3/2}} \nonumber\\
&&-\frac{W_1}{r^2_1}\biggl\{\frac{y}{r^2_1}[(x+\mu)\dot{x}+y\dot{y}]+\dot{y}+n(x+\mu)\biggr\},\end{eqnarray*}

\begin{eqnarray}
&&\Omega=\frac{n^2(x^2+y^2)}{2}+\frac{(1-\mu)q_1}{r_1}+\frac{\mu}{r_2}\nonumber\\&&+\frac{\mu
 A_2}{2r_2^3}+\frac{M_b}{\left(r^2+T^2\right)^{1/2}}\nonumber\\&&+W_1\left[ \frac{(x+\mu)\dot x + y\dot y}{2r_1^2}-n \arctan\left( \frac{y}{x+\mu}\right)\right]\nonumber\\&&\label{eq:OmegaFF}
 \end{eqnarray}
$W_1=\frac{(1-\mu)(1-q_1)}{c_d}$, $q_1=1-\frac{F_p}{F_g}$ is a mass reduction factor expressed in terms of the particle radius $\mathbf{a}$, density $\rho$ radiation pressure efficiency factor $\chi$ (in C.G.S. system):\( q_1=1-\frac{5.6\times{10^{-5}}}{\mathbf{a}\rho}\chi
\). The energy integral of the problem is given by  $C=2\Omega-{\dot{x}}^2-{\dot{y}}^2$, where the quantity $C$ is the Jacobi's  constant. The zero velocity curves $C=2\Omega(x,y)\label{eq:C}$ are presented in various frames of  figure (~\ref{fig:plotzvc})  for the entire range of parameters $q_1, A_2, M_b$ details of the frames are given in the table (~\ref{Tab:1}). We have seen that  in frame $B(q_1=.5)$ and $C(q_1=1)$ there are closed curves around the $L_{4(5)}$ so they are stable but the stability range reduced[see frame $B(q_1=.5)$ in figure (~\ref{fig:plotzvc})] due to P-R effect. In frame $A(q_1=0)$ the closed curves around $L_{4(5)}$ disappeared so they are unstable due to P-R effect. When the  P-R effect is absent as in frames $D(q_1=1,A_2=0.02 I-M_b=0.25,II-M_b=0.5, III-M_b=0.75 )$ the effect of  oblateness and mass of the belt is presented, we have seen that the position of  ovals around the  $L_{4(5)}$ are different but still  $L_{4(5)}$ are stable.
\begin{table}
\begin{center}
\caption[]{Zero velocity curves   when  $ T= 0.01, \mu= 0.025 $}\label{Tab:1}
 \begin{tabular}{|c|l|c|c|c|l|}
  \hline\noalign{\smallskip}
Frame &  I($A_2=0.00 $ ) & II( $A_2=0.02 $) & III($A_2=0.04 $)  &   About the stability of      \\
when $M_b=0.2$   & oval  around $L_{4(5)}$ &   oval  around $L_{4(5)}$   &  oval  around $L_{4(5)}$  & $L_{4(5)}$ \\\hline\noalign{\smallskip} 
A ($q_1=0$ ) & No   &   No   &  No & Unstable  \\
B ($q_1=0.5$) & very small  &  very small  &  very small  & Stability reduced              \\
C ($q_1=1$ )&  yes  & yes       & yes & Stable                \\\hline\noalign{\smallskip}
Frame & & &  & Stability affected by \\
$q_1=1 A_2=0.02$ & I($M_b=0.25$) & II($M_b=0.5 $)& I($M_b=0.75$)& mass of the belt \\\hline\noalign{\smallskip}
D   &  low effect of belt   &  midium effect of belt & very high effect of belt  &    yes           \\
  \noalign{\smallskip}\hline
\end{tabular}
\end{center}
\end{table}

The  position equilibrium points  are given by putting  $\Omega_x=\Omega_y=0$ i.e.,  
\begin{eqnarray}
&& n^{2}x-\frac{(1-\mu)q_1(x+\mu)}{r^3_1}-\frac{\mu(x+\mu-1)}{r^3_2}\nonumber\\&&-\frac{3}{2}\frac{\mu{A_2}(x+\mu-1)}{r^5_2}-\frac{M_b x}{\left(r^2+T^2\right)^{3/2}} \nonumber\\
&&+\frac{W_1 ny}{r^2_1}=0\label{eq:eq1pts},\\
&&n^{2}y
-\frac{(1-\mu)q_{1}{y}}{r^3_1}-\frac{\mu{y}}{r^3_2}-\frac{3}{2}\frac{\mu{A_2}y}{r^5_2}\nonumber\\&&-\frac{M_b y}{\left(r^2+T^2\right)^{3/2}} -\frac{W_1 n(x+\mu)}{r^2_1}=0\label{eq:eq2pts}\end{eqnarray}
when  $(W_1\neq0 )$,  from equations  (~\ref{eq:eq1pts}, ~\ref{eq:eq2pts}) we obtained:
\begin{eqnarray}
r_1&&=q_1^{1/3}\left[1-\frac{nW_1}{6(1-\mu)y}-\frac{A_2}{2}\right.\nonumber\\&&\left.+\frac{(1-2r_c)M_b \left(1-\frac{3\mu A_2}{2(1-\mu)}\right)}{3\left(r_c^2+T^2\right)^{3/2}}\right]\label{eq:r1},\\ 
r_2&&=1+\frac{\mu(1-2r_c)M_b}{3\left(r_c^2+T^2\right)^{3/2}}+\frac{nW_1}{3\mu y} \label{eq:r2}
\end{eqnarray}
From above, we obtained:
\begin{eqnarray} &&x=-\mu\pm\left[\left(\frac{q_1}{n^2}\right)^{2/3}\left[1+\frac{nW_1}{2(1-\mu)y}+\frac{3A_2}{2}\right.\right.\nonumber\\&&\left.\left.-\frac{(1-2r_c)M_b \left(1-\frac{3\mu A_2}{2(1-\mu)}\right)}{\left(r_c^2+T^2\right)^{3/2}}\right]^{-2/3}-y^2\right]^{1/2} \label{eq:x1lag}\\
&&x=1-\mu\nonumber\\&&\pm\Bigl[\left[1-\frac{nW_1}{\mu y}(1-\frac{5}{2}A_2)-\frac{\mu(1-2r_0)M_b}{\left(r_c^2+T^2\right)^{3/2}}\right]^{-2/3}-y^2\Bigr]^{1/2}\label{eq:x2lag}
\end{eqnarray}
 From  equations (~\ref{eq:eq1pts}, ~\ref{eq:eq2pts}) the value of $y$ is allwayes  positive, hence the  the equilibrium points are no-longer collinear with the primaries. 
The triangular equilibrium points are given by putting  $\Omega_x=\Omega_y=0$, $y\neq{0}$,   then from equations (~\ref{eq:Omegax}) and (~\ref{eq:Omegay}) we obtained the triangular equilibrium points as:%
\begin{eqnarray}
&&x=-\mu+\frac{q_1^{2/3}}{2}(1-A_2)-\frac{nW_1\left[\mu q_1^{2/3}-2(1-\mu)\right]}{6\mu(1-\mu)y_0}\nonumber\\&&+ \frac{(1-2r_c) M_b\left[\left\{1-\frac{3\mu A_2}{(1-\mu)}\right\}q_1^{2/3}-1\right]}{3\left(r_c^2+T^2\right)^{3/2}}\label{eq:xl4}\end{eqnarray}
\begin{eqnarray}
 &&y=\pm \frac{q_1^{2/3}}{2} \left[4-q_1^{2/3}+2\left(q_1^{2/3}-2\right)A_2\right.\nonumber\\&&\left.-\frac{2nW_1\left(q_1^{2/3}-2\right)}{ 3\mu(1-\mu)y_0}\right.\nonumber\\&&\left.-\frac{4(2r_c-1)M_b\left[\left\{\left(q_1^{2/3}-3\right)-\frac{3\mu A_2\left(q_1^{2/3}-3\right)}{2(1-\mu)}\right\}\right]}{3\left(r_c^2+T^2\right)^{3/2}}\right]^{1/2}\label{eq:yl4} 
 \end{eqnarray}
All these results are similar with \cite{Szebehely1967}, \cite{Ragosetal1995},   \cite{Kushvah2008Ap&SS, Kushvah2008Ap&SS.315} and others. 
\section{Lyapunov stability}
\label{sec:lstb}
In this section we will examine the P-R effect on the  linear stability conditions. In mathematics, the notion of Lyapunov s stability occurs in the study of dynamical systems. In simple terms, if all solutions of the dynamical system
that start out near an equ equilibrium m point $L_i$ stay near $L_i$ forever, then $L_i$ is
Lyapunov's table.  Let the position of any equilibrium point is $(x*,y*)$  the taking  $x=x*+\alpha$,\,  $y=y*+\beta$, where  $\alpha=\xi e^{\lambda{t}}$,\ $\beta=\eta e^{\lambda{t}}$ are the small displacements  $\xi,\eta$,\  $\lambda$ are parameters,   then the  equations of perturbed motion corresponding to the system of equations (~\ref{eq:Omegax}), (~\ref{eq:Omegay}) are as follows:
\begin{align}
\ddot{\alpha}-2n\dot{\beta} &= {\alpha}{\Omega^*_{xx}}+{\beta}{\Omega^*_{xy}}+\dot{\alpha}{\Omega^*_{x\dot{x}}}+{\dot{\beta}}{\Omega^*_{x\dot{y}}} \\
\ddot{\beta}+2n\dot{\alpha}&= {\alpha}{\Omega^*_{yx}}+{\beta}{\Omega^*_{yy}}+\dot{\alpha}{\Omega^*_{y\dot{x}}}+\dot{\beta}{\Omega^*_{y\dot{y}}}
\end{align}
where superfix $*$ is corresponding to the equilibrium points.
 \begin{align}(\lambda^2-\lambda{\Omega^*_{x\dot{x}}}-{\Omega^*_{xx}})\xi
+[-(2n+{\Omega^*_{x\dot{y}}})\lambda-{\Omega^*_{xy}}]\eta&=0\label{eq:lambda_x}\\
[(2n-{\Omega^*_{y\dot{x}}})\lambda-{\Omega^*_{yx}}]\xi
+(\lambda^2-\lambda{\Omega^*_{y\dot{y}}}-{\Omega^*_{yy}})\eta&=0\label{eq:lambda_y}
\end{align}
this system has singular solution if,
\[
\begin{vmatrix}
\lambda^2-\lambda{\Omega^*_{x\dot{x}}}-\Omega^*_{xx}& -(2n+{\Omega^*_{x\dot{y}}})\lambda-\Omega^*_{xy} \\(2n-{\Omega^*_{y\dot{x}}})\lambda-\Omega^*_{yx}& \lambda^2-\lambda{\Omega^*_{y\dot{y}}}-\Omega^*_{yy}\\
\end{vmatrix}
=0
\]
\begin{equation}
\Rightarrow \quad \lambda^4+a\lambda^3+b\lambda^2+c\lambda+d=0 \label{eq:cheq}
\end{equation}
at the equilibrium points:
\begin{eqnarray*}
&&a=3\frac{W_1}{r^2_{1*}},\\
 &&b=2n^2-f_*-\frac{3\mu{A_2}}{{r^5_2}_*}+\frac{3M_b T^2}{\left(r^2_*+T^2\right)^{5/2}}+\frac{2W^2_1}{{r^4_1}_*}\\
&&c=-a(1+e),\\ && e=\frac{\mu}{{r^5_2}_*}A_2+\frac{\mu}{{{r^2_1}_*}{{r^5_2}_*}}\zabs{y^2_*}\\&&+\frac{3M_b \left(\frac{\mu^2 y^2_*}{r^2_*}-T^2\right)}{\left(r^2_*+T^2\right)^{5/2}}
 \end{eqnarray*}
\begin{eqnarray*}&&d=(n^2-f_*)\biggl[n^2+2f_*-\frac{3\mu{A_2}}{{r^5_2}_*}+\frac{3M_b T^2}{\left(r^2_*+T^2\right)^{5/2}}\biggr]\nonumber\\&&+9\mu\amc{y^2_*}\left[\frac{q_1}{r^5_{1*}r^5_{2*}}\right.\\&&\left.+\frac{3M_b}{\left(r^2_*+T^2\right)^{5/2}}\left\{\frac{\mu q_1}{{r_1^5}_*}+\frac{(1-\mu)\zabs}{r_{2*}^5}\right\}\right]\nonumber\\&&
-\frac{6\mu nW_1y_*}{{{r^4_1}_*}}\left\{\frac{(x_*+\mu)(x_*+\mu-1)+y^2_*}{{r^5_2}_*}\right.\\&&\left.+\frac{3M_b \left[x_*(x_*+\mu)+y^2_*\right]}{\left(r^2_*+T^2\right)^{5/2}}\right\} \end{eqnarray*}
where \(f_*=\frac{(1-\mu)q_1}{r^3_{1_*}}+\frac{\mu}{r^3_{2*}}\left(1+\frac{3}{2}\frac{A_2}{r^2_{2*}}\right)+\frac{3M_b }{\left(r^2_*+T^2\right)^{5/2}}\). The points $L_1,L_2,L_3$  no longer lie along the line joining the primaries, since the condition is not satisfied for them, so taking $y\rightarrow 0,\frac{W_1}{y}\rightarrow 0$ because $y>>W_1, x>>W_1$, from (~\ref{eq:r1}) we have  $r_1\approx \left[\frac{q_1}{n^2}\right]^{1/3}$. Since   $\mathbf{f_*}>1$ and  characteristic equation (~\ref{eq:cheq})  positive root for   collinear points so they are unstable.
 Using Ferrari's theorem the roots of characteristic equation (~\ref{eq:cheq}) are given by:
\begin{equation}
\lambda_i=-\frac{(a+A)}{4}\pm\sqrt{\left(\frac{a+A}{4}\right)^2-B}\label{eq:cheq_withPR}
\end{equation}
where $A=\pm\sqrt{8l-4b+a^2}$, and \( B=\left( \frac{b}{2}+\alpha_1a^2\right)\left(1\pm\sqrt{1+8\alpha_1}\right)\mp\frac{1+e}{\sqrt{1+8\alpha_1}}\), $\alpha_1=\frac{(1+e)(1+e^2-b)+d}{2(b^2-4d)}>0 $ from this the characteristic roots are given by:
\begin{eqnarray}
\lambda_{1,2}&=&-\frac{a\left(1+\sqrt{1+8\alpha_1}\right)}{4}\nonumber\\&&\pm\sqrt{\frac{a^2\left(1+\sqrt{1+8\alpha_1}\right)}{16}-B_1}\label{eq:lambda12}\\
\lambda_{3,4}&=&-\frac{a\left(1-\sqrt{1+8\alpha_1}\right)}{4}\nonumber\\&&\pm\sqrt{\frac{a^2\left(1-\sqrt{1+8\alpha_1}\right)}{16}-B_2}\label{eq:lambda34}
\end{eqnarray}
where
\[ B_1=\left( \frac{b}{2}+\alpha_1a^2\right)\left(1+\sqrt{1+8\alpha_1}\right)-\frac{1+e}{\sqrt{1+8\alpha_1}},\]
\[ B_2=\left( \frac{b}{2}+\alpha_1a^2\right)\left(1-\sqrt{1+8\alpha_1}\right)+\frac{1+e}{\sqrt{1+8\alpha_1}}\]
 from  equations(~\ref{eq:lambda12},~\ref{eq:lambda34}), we found that at least one of the roots $\lambda_i(i=1,2,3,4)$  have a positive real part due to P-R effect thus  the triangular  equilibrium points are unstable in the sense of Lyapunov stability  result are similar to \cite{Chernikov1970}  and  \cite{Kushvah2008Ap&SS.315}.
\section{Conclusion}
We have seen the  collinear points $L_1,L_2,L_3$  no longer lie along the line joining the primaries and they are  unstable in classical case. If   $ q_1=.5, q_1=1 $ there are closed curves around the $L_{4(5)}$ so they are stable but the stability range reduced[see frame $B(q_1=0.5)$ in figure (~\ref{fig:plotzvc})] due to P-R effect.  When the  P-R effect is absent as in frames $D(q_1=1,A_2=0.02 I-M_b=0.25,II-M_b=0.5, III-M_b=0.75 )$ the effect of  oblateness and mass of the belt is presented, we have seen that the positions of  ovals around the  $L_{4(5)}$ are different but still  $L_{4(5)}$ are stable. In frame $A(q_1=0)$ the closed curves around $L_{4(5)}$ disappeared so they are unstable due to P-R effect. 
\label{sec:con}
\bibliographystyle{spbasic} 

\end{document}